\documentclass[a4paper,12pt]{amsart}

\usepackage{enumerate, booktabs}
\usepackage{amsmath, amsfonts, amssymb, amsthm}
\usepackage[all]{xy}
\usepackage[height=22.5cm, width=15.5cm]{geometry}

\numberwithin{equation}{section}

\theoremstyle{plain}
\newtheorem{thm}{Theorem}[section]
\newtheorem{prop}[thm]{Proposition}
\newtheorem{lem}[thm]{Lemma}
\newtheorem{cor}[thm]{Corollary}
\newtheorem*{thm*}{Theorem}

\theoremstyle{definition}
\newtheorem{defn}[thm]{Definition}

\theoremstyle{remark}
\newtheorem*{rem}{Remark}
\newtheorem*{ex}{Example}

\newcommand{\mbb}[1]{\mathbb{#1}}
\newcommand{\wt}[1]{\widetilde{#1}}
\newcommand{\wh}[1]{\widehat{#1}}
\newcommand{\ol}[1]{\overline{#1}}
\newcommand{\lie}[1]{{\mathfrak{#1}}}

\newcommand{\hq}{/\hspace{-0.12cm}/}

\DeclareMathOperator{\Alb}{Alb}

\title{Homogeneous K\"ahler and Hamiltonian manifolds}
\thanks{We gratefully acknowledge that this work was partially supported by an
NSERC Discovery Grant. We would also like to thank Nicolas Dutertre for
valuable discussions on subanalytic geometry.}
\dedicatory{Dedicated to Alan T. Huckleberry}

\date{June 2, 2010}


\author{Bruce Gilligan}
\address{Department of Mathematics and Statistics, University of Regina,
Regina, Canada S4S 0A2}
\email{gilligan@math.uregina.ca}

\author{Christian Miebach}
\address{LATP-UMR(CNRS) 6632, CMI-Universit\'e d'Aix-Marseille I, 39, rue
Joliot-Curie, F-13453 Marseille Cedex 13, France}
\email{miebach@cmi.univ-mrs.fr}

\author{Karl Oeljeklaus}
\address{LATP-UMR(CNRS) 6632, CMI-Universit\'e d'Aix-Marseille I, 39, rue
Joliot-Curie, F-13453 Marseille Cedex 13, France}
\email{karloelj@cmi.univ-mrs.fr}

\begin{document}

\maketitle

\begin{abstract}
We consider actions of reductive complex Lie groups $G=K^\mbb{C}$ on
K\"ahler manifolds $X$ such that the $K$--action is Hamiltonian and prove then
that the closures of the $G$--orbits are complex-analytic in $X$. This is
used to characterize reductive homogeneous K\"ahler manifolds in terms of their
isotropy subgroups. Moreover we show that such manifolds admit $K$--moment maps
if and only if their isotropy groups are algebraic.
\end{abstract}

\section{Introduction}
A reductive complex Lie group $G$ is a complex Lie group admitting a compact
real form $K$, i.\,e.\ $G=K^\mbb{C}$. Equivalently a finite covering of $G$ is
of the form $S\times Z =S \times(\mbb{C}^*)^k$, where $S$ is a semisimple
complex Lie group. It is well known that every complex reductive Lie group
admits a unique structure as a linear algebraic group. Holomorphic or algebraic
actions of reductive Lie groups appear frequently in complex and algebraic
geometry and interesting connections arise between the structure of the orbits
of such groups and the isotropy subgroups of the orbits.

A result of this type was proved independently by Matshushima \cite{Mat} and
Onishchik \cite{On}. They consider $G$ a complex reductive Lie group and $H$ a
closed complex subgroup of $G$ and show that $G/H$ is Stein if and only if $H$
is a reductive subgroup of $G$. In \cite{BaOt} Barth and Otte prove that the
holomorphic separability of the homogeneous space $G/H$ implies $H$ is an
algebraic subgroup of the reductive group $G$.

In the case of semisimple actions, it is known that K\"ahler is equivalent to
algebraic in the sense that $S/H$ is K\"ahler if and only if $H$ is an
algebraic subgroup of the complex semisimple Lie group $S$, see \cite{Bert} and
\cite{BerO}. The simple example of an elliptic curve $\mbb{C}^*/\mbb{Z}$ shows
that this result does not hold in the reductive case. Instead homogeneous
K\"ahler manifolds $X=G/H$ with $G=S\times(\mbb{C}^*)^k$ reductive are
characterized by the two conditions $S\cap H$ is algebraic and $SH\subset G$ is
closed, as we shall prove. If, in addition to the existence of a K\"ahler form,
there exists a $K$--moment map on $X$, then $X$ is called a  {\it Hamiltonian}
$G$--manifold. Huckleberry has conjectured that the isotropy groups in
a Hamiltonian $G$--manifold are algebraic. In the present paper we prove that
this is indeed the case.

The moment map plays a decisive role in our proof which depends in an essential
way on the work of Heinzner-Migliorini-Polito \cite{HeMiPo}. In the third
section of their paper they investigate the closure of certain orbits and prove
the following:  suppose $T$ is an algebraic torus acting holomorphically on a
complex space $X$ such that the semistable quotient $\pi\colon X\to X\hq T$
exists. Let $A$ be a subanalytic set in $X$ such that $\pi|_A\colon A\to X\hq
T$ is proper. Then $T\cdot A$ is subanalytic in $X$. We use the moment map in
order to ensure the existence of the semistable $T$--quotient locally. This is
sufficient to show that the $G$--orbits are locally subanalytic and hence
locally closed in the Hamiltonian $G$--manifold $X$. Moreover, we deduce
from this fact that the closure of any $G$--orbit is complex-analytic in $X$.
This generalizes previous work of~\cite{Som} and~\cite{Fu} to non-compact
K\"ahler manifolds.

Our work was partially motivated by~\cite{Marg}, where Margulis constructed
discrete subgroups $\Gamma$ of ${\rm{SL}}(2,\mbb{R})\ltimes\mbb{R}^3$ which are
free groups generated by two elements. These groups $\Gamma$ can be divided
into two non-empty classes depending on whether the induced action of $\Gamma$
on $\mbb{R}^3$ is properly discontinuous or not. The associated homogeneous
complex manifolds $({\rm{SL}}(2,\mbb{C})\ltimes\mbb{C}^3)/\Gamma$ are not
K\"ahler in the ``non properly discontinuous'' case by Corollary
\ref{Cor:locallyclosed}. It seems to be a difficult problem to decide the
K\"ahler question for these quotients in the ``properly discontinuous'' case.

The paper is organized as follows. In section $2$, the definitions of
$K$--moment maps and Hamiltonian actions are recalled. Furthermore two lemmata
are proved for later use. The main result of section $3$ is the analyticity of
orbit closures. Since for $G=S$ semisimple there always is a moment map, any
semisimple Lie group action on a K\"ahler manifold has locally closed orbits.

In section $4$ we prove that the reductive homogeneous manifold $X=G/H$ is
Hamiltonian if and only if $H$ is an algebraic subgroup of $G$ and use this to
give a new proof of the main results in~\cite{BerO} and~\cite{Bert}. Finally,
in the last section this result is used in order to prove our characterization
of those closed complex subgroups $H\subset G$ such that $X=G/H$ admits a
K\"ahler form.

\section{Hamiltonian $G$--manifolds}

Let $G=K^\mbb{C}$ be a complex reductive  Lie group with maximal compact
subgroup $K$. Let $X$ be a complex manifold endowed with a holomorphic
$G$--action.

We denote the Lie algebra of $K$ by $\lie{k}$. The group $K$ acts via the
coadjoint representation on the dual $\lie{k}^*$. In the following equivariance
of a map with values in $\lie{k}^*$ is always meant with respect to the
coadjoint action. If $\xi\in\lie{k}$, we write $\xi_X$ for the holomorphic
vector field on $X$ whose flow is given by $(t,x)\mapsto\exp(t\xi)\cdot x$. If
$\omega$ is a $K$--invariant K\"ahler form on $X$, then the contracted form
$\iota_{\xi_X}\omega$ is closed for every $\xi\in\lie{k}$. By definition, a
$K$--equivariant smooth map $\mu\colon X\to\lie{k}^*$ is a moment map for the
$K$--action on $X$ if for each $\xi\in\lie{k}$ the smooth function $\mu^\xi\in
\mathcal{C}^\infty(X)$, $\mu^\xi(x):=\mu(x)\xi$, verifies $d\mu^\xi=
\iota_{\xi_X}\omega$. The $K$--action on $X$ is called {\emph{Hamiltonian}} if
an equivariant moment map $\mu\colon X\to\lie{k}^*$ exists. Note that, if $\mu$
is a moment map and if $\lambda\in\lie{k}^*$ is a $K$--fixed point, then
$\mu+\lambda$ is another moment map on $X$.

\begin{defn}\label{Defn:Hamiltonian}
We say that $X$ is a {\emph{Hamiltonian $G$--manifold}} if $X$ admits a
$K$--invariant K\"ahler form such that the $K$--action on $X$ is Hamiltonian
with equivariant moment map $\mu\colon X\to\lie{k}^*$.
\end{defn}

\begin{rem}
If $G$ is semisimple, then every K\"ahler manifold $X$ on which $G$ acts
holomorphically is a Hamiltonian $G$--manifold which can be seen as follows.
Let $dk$ be the normalized Haar measure of $K$. If $\omega$ is any K\"ahler
form on $X$, then $\wh{\omega}:=\int_K k^*\omega dk$ is a $K$--invariant
K\"ahler form on $X$. Since $K$ is semisimple, there exists a unique
equivariant moment map $\mu\colon X\to\lie{k}^*$ by Theorem~26.1
in~\cite{GS}.
\end{rem}

In this paper we will often use the following.

\begin{lem}\label{Lem:restriction}
Let $X$ be a Hamiltonian $G$--manifold and let $\wt{G}$ be a complex reductive 
subgroup of $G$. Then every $\wt{G}$--stable complex submanifold $\wt{X}$ of
$X$ is a Hamiltonian $\wt{G}$--manifold.
\end{lem}

\begin{proof}
We may assume without loss of generality that $\wt{G}=\wt{K}^\mbb{C}$ for some
compact subgroup $\wt{K}\subset K$. Composing the moment map $\mu\colon X\to
\lie{k}^*$ with the orthogonal projection onto $\wt{\lie{k}}^*$ we obtain a
moment map for the $\wt{K}$--action on $X$. Restricting this map to the
K\"ahler manifold $\wt{X}$ we see that the $\wt{K}$--action on $\wt{X}$ is
Hamiltonian.
\end{proof}

\begin{ex}
Let $G\to{\rm{GL}}(V)$ be a holomorphic representation of the complex reductive 
group $G$ on a finite dimensional complex vector space $V$. Then each
$G$--stable complex submanifold of $V$ or of $\mbb{P}(V)$ is a Hamiltonian
$G$--manifold. In particular, if $H$ is an algebraic subgroup of $G$, then the
homogeneous space $G/H$ is a quasi-projective variety (see e.\,g.\ Theorem~5.1
in~\cite{Bo}) and hence a Hamiltonian $G$--manifold.
\end{ex}

For later use we note the following

\begin{lem}\label{Lem:cover}
Let $(X,\omega)$ be a Hamiltonian $G$-manifold
with $\mu:X \to \mathfrak k^*$ its moment map and let $p\colon\wt{X}\to X$ be
a topological covering. If the $G$--action lifts to $\wt{X}$, then
$(\wt{X},p^*\omega)$ is a Hamiltonian $G$--manifold with moment map $p^*\mu$.
\end{lem}

\begin{proof}
We equip $\wt{X}$ with the unique complex structure such that $p$ is locally
biholomorphic. If the $G$--action lifts to $\wt{X}$, then $G$ acts
holomorphically on $\wt{X}$ and $p$ is $G$--equivariant. Consequently,
$p^*\omega$ is a $K$--invariant K\"ahler form on $\wt{X}$.

For $\xi\in\lie{k}$ let $\xi_{\wt{X}}$ and $\xi_X$ be the induced vector fields
on $\wt{X}$ and $X$, respectively. Since $p$ is equivariant, we have $p_*
\xi_{\wt{X}}=\xi_X$. Hence, we obtain
\begin{equation*}
d(p^*\mu)^\xi=dp^*\mu^\xi=p^*d\mu^\xi=p^*\iota_{\xi_X}\omega=
\iota_{\xi_{\wt{X}}}p^*\omega,
\end{equation*}
which shows that $p^*\mu$ is an equivariant moment map for the $K$--action on
$\wt{X}$.
\end{proof}

\section{Local closedness of $G$--orbits}

Let $X$ be a Hamiltonian $G$--manifold where $G=K^\mbb{C}$ is a complex
reductive group. Suppose that $X$ is $G$--connected, i.\,e.\ that $X/G$ is
connected. In particular, $X$ has only finitely many connected components since
this is true for $G$. We want to show that the topological closure of every
$G$--orbit is complex-analytic in $X$.

We fix a maximal torus $T_0$ in $K$. Then $T:=T_0^\mbb{C}$ is a maximal
algebraic torus in $G$ and the moment map $\mu\colon X\to\lie{k}^*$ induces by
restriction a moment map $\mu_T\colon X\to\lie{t}_0^*$ for the $T_0$--action on
$X$. Since $\lie{t}_0$ is Abelian, for every $\lambda\in\lie{t}_0^*$ the
shifted map $\mu_T+\lambda$ is again a moment map for $T_0$. Consequently,
every $x\in X$ lies in the zero fiber of some moment map for the $T_0$--action
on $X$ which has the following consequences (see~\cite{HeLo}).

\begin{thm}\label{Thm:SliceQuotient}
Let $G$ be a complex reductive  group and $X$ be a $G$--connected Hamiltonian
$G$--manifold.
\begin{enumerate}
\item Every isotropy group $T_x$ is complex reductive and hence the connected
component of the identity $(T_x)^0$ is a subtorus of $T$.
\item For every $x\in X$ there exists a complex submanifold $S$ of $X$ which
contains $x$ such that the map $T\times_{T_x}S\to T\cdot S$, $[t,y]\mapsto
t\cdot y$, is biholomorphic onto its open image.
\item For $\lambda\in\lie{t}_0^*$ we define $\mathcal{S}_\lambda:=\bigl\{x\in
X;\ \ol{T\cdot x}\cap\mu_T^{-1}(\lambda)\not=\emptyset\bigr\}$. Then
$\mathcal{S}_\lambda$ is a $T$--stable open subset of $X$ such that the
semistable quotient (see~\cite{HeMiPo}) $\mathcal{S}_\lambda\to
\mathcal{S}_\lambda\hq T$ exists. Moreover, the inclusion $\mu_T^{-1}(\lambda)
\hookrightarrow\mathcal{S}_\lambda$ induces a homeomorphism $\mu_T^{-1}(
\lambda)/T_0\cong\mathcal{S}_\lambda\hq T$.
\end{enumerate}
\end{thm}

\begin{rem}
Properties $(1)$ and $(2)$ imply that if the $T$--action on $X$ is known to be
almost free, then it is locally proper.
\end{rem}

For the following we have to review the definition of subanalytic sets. For
more details we refer the reader to~\cite{BiMil} and to~\cite{Hi}.

Let $M$ be a real analytic manifold. A subset $A\subset M$ is called
{\emph{semianalytic}} if every point in $M$ has an open neighborhood $\Omega$
such that $A\cap\Omega=\bigcup_{k=1}^r\bigcap_{l=1}^sA_{kl}$, where every
$A_{kl}$ is either of the form $\{f_{kl}=0\}$ or $\{f_{kl}>0\}$ for
$f_{kl}\in\mathcal{C}^\omega(\Omega)$. A subset $A\subset M$ is called
{\emph{subanalytic}} if every element of $M$ admits an open neighborhood
$\Omega$ such that $A\cap\Omega$ is the image of a semianalytic set under a
proper real analytic map. We note that finite intersections and finite unions
as well as topological closures of subanalytic sets are subanalytic. Finally we
call a set $A\subset M$ {\emph{locally subanalytic}} if there are open sets
$U_1,\dotsc,U_k\subset M$ such that $A\subset U_1\cup\dotsb\cup U_k$ and such
that $A\cap U_j$ is subanalytic in $U_j$ for every $j$. For later use we cite
the following theorem of Hironaka (\cite{Hi}).

\begin{thm}\label{Thm:Hironaka}
Let $\Phi\colon M\to N$ be a real analytic map between real analytic manifolds
and let $A\subset M$ be subanalytic. If $\Phi|_{\ol{A}}\colon\ol{A}\to N$ is
proper, then $\Phi(A)$ is subanalytic in $N$.
\end{thm}

It is shown in~\cite{HeMiPo} that, if the semistable quotient $X\to X\hq T$
exists globally, then the semistable quotient $X\to X\hq G$ exists. The first
step in the proof of this theorem consists in showing that the existence of
$X\hq T$ implies that the $G$--orbits are subanalytic and thus locally closed
in $X$. In our situation the semistable quotient of $X$ with respect to $T$
exists only locally (in the sense of Theorem~\ref{Thm:SliceQuotient}(3)). As we
will see this implies that the $G$--orbits in $X$ are locally subanalytic which
is sufficient for them to be locally closed.

The following lemma is the essential ingredient in the proof of this statement.

\begin{lem}\label{Lem:T}
Let $A\subset X$ be a compact subanalytic set. Then $T\cdot A$ is locally
subanalytic in $X$.
\end{lem}

\begin{proof}
Since $A$ is compact, we have $A\subset\bigcup_{k=1}^n\mathcal{S}_{\lambda_k}$.
For every $k=1,\dotsc,n$ let $U_k$ be an open subanalytic subset of
$\mathcal{S}_{\lambda_k}$ such that $\ol{U}_k\subset\mathcal{S}_{\lambda_k}$ is
compact and such that $A\subset\bigcup_{k=1}^nU_k$. Consequently, for every $k$
the intersection $A\cap\ol{U}_k$ is a compact subanalytic subset of
$\mathcal{S}_{\lambda_k}$. Since for each $k$ the semistable quotient
$\mathcal{S}_{\lambda_k}\to\mathcal{S}_{\lambda_k}\hq T$ exists, we conclude
from the proposition in Section~3 of~\cite{HeMiPo} that
$T\cdot(A\cap\ol{U}_k)=(T\cdot A)\cap(T\cdot\ol{U}_k)$ is subanalytic in
$\mathcal{S}_{\lambda_k}$. It follows that for every $k$ the intersection
$(T\cdot A)\cap(T\cdot U_k)$ is subanalytic in the open set $T\cdot U_k\subset
X$. Since we have $T\cdot A=T\cdot\bigl(\bigcup_{k=1}^nA\cap U_k\bigr)=
\bigcup_{k=1}^n\bigl((T\cdot A)\cap(T\cdot U_k)\bigr)$, we conclude that
$T\cdot A$ is locally subanalytic in $X$.
\end{proof}

\begin{lem}\label{Lem:K}
Let $A\subset X$ be (locally) subanalytic. Then $K\cdot A$ is (locally)
subanalytic in $X$.
\end{lem}

\begin{proof}
Since $K$ is compact, the real analytic map $\Phi\colon K\times X\to X$,
$(k,x)\mapsto k\cdot x$, is proper: For every compact subset $C\subset X$ the
inverse image $\Phi^{-1}(C)$ is closed and contained in $K\times(K\cdot C)$,
hence compact. We conclude that the restriction of $\Phi$ to $K\times\ol{A}$ is
proper. Therefore Hironaka's theorem~\ref{Thm:Hironaka}, \cite{Hi} implies that
$\Phi(K\times A)=K\cdot A$ is subanalytic.

If $A$ is locally subanalytic, then $A$ is covered by relatively compact
subanalytic open sets $U$ such that $A\cap\ol{U}$ is subanalytic. Then it
follows as above that $K\cdot(A\cap\ol{U})$ is subanalytic, and consequently
$K\cdot A$ is locally subanalytic.
\end{proof}

Now we are in a position to prove the main result of this section.

\begin{thm}\label{Thm:locallyclosed}
Suppose $X$ is a $G$--connected Hamiltonian $G$--manifold, where $G$ is a
complex reductive  group. Then
\begin{enumerate}
\item every $G$--orbit is locally subanalytic and in particular locally closed
in $X$,
\item the boundary of every $G$--orbit contains only $G$--orbits of strictly
smaller dimension, and
\item the closure of every $G$--orbit is complex-analytic in $X$.
\end{enumerate}
\end{thm}

\begin{proof}
For every $x\in X$ the orbit $K\cdot x$ is a compact real analytic submanifold
of $X$. By Lemma~\ref{Lem:T} the set $T\cdot(K\cdot x)$ is locally subanalytic
in $X$. Thus Lemma~\ref{Lem:K} implies that $K\cdot\bigl(T\cdot(K\cdot x)
\bigr)$ is locally subanalytic as well. Because of $G=KTK$ every $G$--orbit is
locally subanalytic.

In order to see that the $G$--orbits are locally closed, we take
$U_1\cup\dotsb\cup U_k$ to be an open covering of $G\cdot x$ such that for
every $j$ the intersection $(G\cdot x)\cap U_j$ is subanalytic in $U_j$. Since
the boundary of a subanalytic set is again subanalytic and of strictly smaller
dimension, we see that $(G\cdot x)\cap U_j$ contains an interior point of its
closure in $U_j$. Moving this point with the $G$--action it follows that
$(G\cdot x)\cap U_j$ is open in its closure in $U_j$. Consequently, $G\cdot x$
is locally closed.

For the second claim it is sufficient to note that the dimension of an orbit
$G\cdot x$ can be checked in the intersection with an open set $U$ such that
$(G\cdot x)\cap U$ is subanalytic in $U$. More precisely, let $x,y\in X$ such
that $G\cdot y\subset\ol{G\cdot x}$. Since $\{x,y\}$ is compact subanalytic,
there are finitely many open sets $U_1,\dotsc,U_k$ such that $(G\cdot x)\cup
(G\cdot y)\subset U_1\cup\dotsb\cup U_k$ and such that $\bigl((G\cdot x)\cup
(G\cdot y)\bigr)\cap U_j$ is subanalytic in $U_j$ for every $j$. Suppose
$y\in U_1$. Then $(G\cdot y)\cap U_1$ lies in the closure of $(G\cdot x)\cap
U_1$ in $U_1$. After possibly shrinking $U_1$ we may assume that $(G\cdot y)
\cap U_1$ is subanalytic in $U_1$ which implies that $(G\cdot x)\cap U_1$ is
also subanalytic in $U_1$. Hence, we obtain $\dim G\cdot y=\dim(G\cdot y)\cap
U_1<\dim(G\cdot x)\cap U_1=\dim G\cdot x$ as was to be shown.

Finally let $x_0 \in X$ and $E:=\{x\in X;\ \dim
G\cdot x < \dim G \cdot {x_0} \}$. The set $E$ is complex-analytic
and its complement $\Omega:= X \setminus E$ is $G$--invariant.
Since the boundary of $G\cdot x_0$ contains only orbits of strictly smaller
dimension by the previous claim, the orbit $G\cdot x_0$ is closed in $\Omega$
and therefore a complex
submanifold of $\Omega$. We will show that $\ol{G\cdot x_0}$ is
complex-analytic in $X$ by applying Bishop's theorem (\cite{Bi}). For this we
must check that every point $x\in E$ has an open neighborhood $U\subset X$ such
that $U\cap(G\cdot x_0)$ has finite volume with respect to some hermitian
metric on $X$. Without loss of generality we may assume that $x\in\ol{G\cdot
x_0}\cap E$ holds. According to what we have already shown we find an open
neighborhood $U\subset X$ of $x$ such that $U\cap\ol{G\cdot x_0}$ is
subanalytic in $U$. After possibly shrinking $U$ we may assume that $U$ is
biholomorphic to the unit ball in $\mbb{C}^n$. It is known (see the remark
following Proposition~1.4 in~\cite{KuRa}) that the $2k$--dimensional Hausdorff
volume (where $k:=\dim_\mbb{C}G\cdot x_0$) of $U\cap(G\cdot x_0)$ is finite.
Since $U \cap(G\cdot x_0)$ is also an immersed submanifold of $U$, the
$2k$--dimensional Hausdorff volume coincides with the geometric volume
associated with the standard hermitian metric on $\mbb{C}^n$ (see page~48
in~\cite{Si}). This observation allows us to deduce from Bishop's theorem that
$\ol{G\cdot x_0}$ is complex-analytic in $X$.
\end{proof}

\begin{rem}
We restate the following fact which is shown in the third part of the proof and
might be of independent interest: Let $E \subset \mathbb B_n$ be a
complex-analytic subset. Suppose that $A \subset \mathbb B_n\setminus E$ is
complex-analytic and that $A \subset \mathbb B_n$ is locally subanalytic and an
injectively immersed complex submanifold. Then the topological closure of $A$ in
$\mathbb B_n$ 
is complex-analytic.
\end{rem}
\begin{rem}
In~\cite{Som} holomorphic actions of complex reductive groups $G$ on compact
K\"ahler manifolds $X$ are considered. Under the additional assumption that the
$G$--action on $X$ is projective it is shown that for every $x\in X$ the
closure $\ol{G\cdot x}$ is complex-analytic in $X$. Sommese's notion of
projectivity of a $G$--action on $X$ is equivalent to the fact that $G$ acts
trivially on the Albanese torus $\Alb(X)$. Hence, by Proposition~1 on page~269
in~\cite{HuWu}, $G$ acts projectively on $X$ if and only if $X$ is a
Hamiltonian $G$--manifold.

In~\cite{Fu} some properties of algebraic group actions are extended
to the more general class $\mathcal{C}$ of compact complex spaces
that are the meromorphic images of compact K\"{a}hler spaces and
it is shown that the orbit closures are complex analytic in this setting.
\end{rem}

From the remark after Definition~\ref{Defn:Hamiltonian} we obtain the
following

\begin{cor}\label{Cor:locallyclosed}
Let $G=S$ be a semisimple complex Lie group acting holomorphically on the
K\"ahler manifold $X$. Then the $S$--orbits are locally closed in $X$.
\end{cor}

\section{Homogeneous Hamiltonian $G$--manifolds}

Let $G=K^\mbb{C}$ be a connected complex reductive group and let $H$ be a
closed complex subgroup of $G$. Suppose that the homogeneous space $X=G/H$
admits a $K$--invariant K\"ahler form $\omega$. We want to show that the
existence of a $K$--equivariant moment map $\mu\colon X\to\lie{k}^*$ implies
that $H$ is an algebraic subgroup of $G$.

\begin{ex}
If $G$ is Abelian, i.\,e.\ if $G=(\mbb{C}^*)^k$, then the fact that $G/H$ is
K\"ahler does not imply that $H$ is algebraic as the example of an elliptic
curve $\mbb{C}^*/\mbb{Z}$ shows. However, if $G/H$ is a Hamiltonian
$G$--manifold, then by Theorem~\ref{Thm:SliceQuotient}(1) the group $H$ is
complex reductive  and hence algebraic.
\end{ex}

\begin{ex}
Suppose that $X=G/H$ is a Hamiltonian $G$--manifold with moment map $\mu\colon
X\to\lie{k}^*$. If $\mu^{-1}(0)\not=\emptyset$, then the semistable quotient
$X\hq G$ exists (and is a point) and thus $X=G/H$ is Stein by~\cite{HeMiPo}.
In this case $H$ is a reductive complex subgroup of $G$  and hence is
algebraic, see \cite{Mat} and \cite{On}.
\end{ex}

We will need the following technical result.

\begin{lem}\label{Lem:twistedproduct}
Let $\Gamma$ be a discrete subgroup of $G$ normalizing $H$ such that
$X=(G/H)/\Gamma$ is a Hamiltonian $G$--manifold with moment map $\mu$. Suppose
that $\Gamma$ acts by holomorphic transformations on a complex manifold $Y$ and
that $Y$ admits a $\Gamma$--invariant K\"ahler form $\omega_Y$. Recall that the
twisted product $(G/H)\times_\Gamma Y$ is by definition the quotient of $(G/H)
\times Y$ by the diagonal $\Gamma$--action $\gamma\cdot(gH,y):=(\gamma\cdot gH,
\gamma\cdot y)$. Then $(G/H)\times_\Gamma Y$ is a Hamiltonian $G$--manifold
with moment map $\wh{\mu}\colon(G/H)\times_\Gamma Y\to\lie{k}^*$,
$\wh{\mu}[gH,y]:=\mu(gH\Gamma)$.
\end{lem}

\begin{proof}
Let $p\colon G/H\to(G/H)/\Gamma$ be the quotient map and let $\omega$ be a
$K$--invariant K\"ahler form on $(G/H)/\Gamma$. Then $p^*\omega$ is a
$K$-- and $\Gamma$--invariant K\"ahler form on $G/H$ and thus
$p^*\omega+\omega_Y$ is a $\Gamma$--invariant K\"ahler form on $(G/H)\times Y$.
Hence, we see that $\wh{Y}:=(G/H)\times_\Gamma Y$ is K\"ahler.

The map $\wh{\mu}$ is well-defined and $K$--equivariant. Let $\xi_{\wh{Y}}$ be
the vector field on $\wh{Y}$ induced by $\xi\in\lie{k}$. Let $U\subset
(G/H)/\Gamma$ be an open set such that the bundle $q\colon\wh{Y}=
(G/H)\times_\Gamma Y\to(G/H)/\Gamma$ is trivial over $U$, i.\,e.\ such that
$q^{-1}(U)\cong U\times Y$. For every $[gH,y]\in U$ the vector
$\xi_{\wh{Y}}[gH,y]$ corresponds to $\bigl(\xi_{(G/H)/\Gamma}
\bigl(p(gH)\bigr),0\bigr)\in T_{p(gH)}(G/H)/\Gamma\oplus T_yY$. Moreover, we
have $d\wh{\mu}^\xi=d\mu^\xi$ in this trivialization. By construction of the
K\"ahler form on $\wh{Y}$ we conclude that $\wh{\mu}$ is a moment map for the
$K$--action on $\wh{Y}=(G/H)\times_\Gamma Y$.
\end{proof}

We will first prove the algebraicity of $H$ under the assumption that $H$ is a
discrete subgroup of $G$. In this case we write $\Gamma$ instead of $H$.

\begin{prop}
Let $G$ be a connected complex reductive  group and let $\Gamma$ be a discrete
subgroup of $G$ such that $X=G/\Gamma$ is a Hamiltonian $G$--manifold. Then
$\Gamma$ is finite.
\end{prop}

\begin{proof}
Let us briefly recall the Jordan decomposition of elements in the affine
algebraic group $G=K^\mbb{C}$ (see Chapter I.4 in~\cite{Bo}). Suppose that $G$
is a subgroup of ${\rm{GL}}(N,\mbb{C})$. An element $\gamma\in G$ is called
semisimple if the matrix representing $\gamma$ is diagonalizable, and unipotent
if the matrix $\gamma-I_N$ is nilpotent. It can be shown that these notions do
not depend on the chosen embedding $G\hookrightarrow{\rm{GL}}(N,\mbb{C})$.
Moreover, every element $\gamma\in G$ has a (unique) Jordan decomposition
$\gamma=\gamma_{\sf s}\gamma_{\sf u}=\gamma_{\sf u}\gamma_{\sf s}$ in $G$,
where $\gamma_{\sf s}$ is semisimple and $\gamma_{\sf u}$ is unipotent.

Suppose there is an element $\gamma\in\Gamma$ with $\gamma_{\sf u}\not=e$.
Then there exists a nilpotent element $\xi\in\lie{g}$ with $\gamma_{\sf u}=
\exp(\xi)$. Since the  group $\exp(\mbb{C}\xi)$ is closed in $G$, the same
holds for the cyclic group $\langle\gamma\rangle:=\{\gamma^m;
\ m\in\mbb{Z}\}\cong\mbb{Z}$. Lemma~\ref{Lem:cover} implies
then that $G/\langle\gamma\rangle$ is a Hamiltonian $G$--manifold. The group
$\langle\gamma\rangle$ acts on $\mbb{C}^*$ by $\gamma^m\cdot z:=e^{im}z$.
Applying Lemma~\ref{Lem:twistedproduct} we conclude that the twisted product
$G\times_{\langle\gamma\rangle}\mbb{C}^*$ is a Hamiltonian $G$--manifold. Since
the $G$--orbits in this twisted product intersect $\mbb{C}^*$ in
$\langle\gamma\rangle$--orbits and since these orbits are dense in the
$S^1$--orbits in $\mbb{C}^*$, we arrive at a contradiction to
Theorem~\ref{Thm:locallyclosed}. Consequently, every $\gamma\in\Gamma$ must be
semisimple.

If $\gamma=\gamma_{\sf s}$, then the Zariski closure of the cyclic group
generated by $\gamma$ is either finite or a complex torus $T\cong(\mbb{C}^*)^l$
for some $l\geq1$. Assume that the latter holds. Then $T/(\Gamma\cap T)$ is a
Hamiltonian $T$--manifold by Lemma~\ref{Lem:restriction}, and consequently
$\Gamma\cap T$ must be finite. Since $\langle\gamma\rangle$ is contained in
$\Gamma\cap T$, it follows that $T$ is finite, a contradiction. Thus every
element of $\Gamma$ is semisimple and generates a finite group. According to
Lemma~2.1 in~\cite{BaOt} the group $\Gamma$ is finite.
\end{proof}

\begin{rem}
If the group $G$ is semisimple, then every holomorphic $G$--manifold which
admits a K\"ahler form is Hamiltonian. Hence, we have given a new proof for
Theorem~3.1 in~\cite{BerO}.
\end{rem}

Now we return to the general case that $H$ is any closed complex subgroup of
$G$ such that $X=G/H$ is a Hamiltonian $G$--manifold. The following theorem
the proof of which can be found in~\cite{BaOt} gives a necessary and sufficient
condition for $H$ to be algebraic.

\begin{thm}
For $h\in H$ let $\mathcal{A}(h)$ denote the Zariski closure of the cyclic
group generated by $h$ in $G$. The group $H$ is algebraic if and only if
$\mathcal{A}(h)$ is contained in $H$ for every $h\in H$.
\end{thm}

Using this result we now prepare the proof of our main theorem in this section.

Let $h\in H$. In order to have better control over the group $\mathcal{A}(h)$
we follow closely an idea which is described on page~107 in~\cite{BaOt}. For
this let $h=h_{\sf s}h_{\sf u}$ be the Jordan decomposition of $h$ in $G$. As
we already noted above, if $h$ is semisimple, then $\mathcal{A}(h)$ is either
finite or isomorphic to $(\mbb{C}^*)^l$. In the first case we have
$\mathcal{A}(h)\subset H$. In the second case, $X=G/H$ is a Hamiltonian
$\mathcal{A}(h)$--manifold which implies that the orbit $\mathcal{A}(h)\cdot
eH\cong\mathcal{A}(h)/\bigl(\mathcal{A}(h)\cap H\bigr)$ is Hamiltonian. Hence,
$\mathcal{A}(h)\cap H$ is algebraic which yields $\mathcal{A}(h)\subset H$.

If $h$ is unipotent, then there exists a simple three dimensional closed
complex subgroup $S$ of $G$ containing $h$ (see~\cite{Jac}). Again $X=G/H$ is a
Hamiltonian $S$--manifold. Hence the orbit $S\cdot eH\cong S/(S\cap H)$ is
Hamiltonian and in particular K\"ahler. We have to show that $S\cap H$ is
algebraic in $S$. Then we have $\mathcal{A}(h)\subset S\cap H\subset H$, as was
to be shown. Algebraicity of $S\cap H$ will be a consequence of the following
lemma for which we give here a direct proof.

\begin{lem}
Let $H$ be a closed complex subgroup of $S={\rm{SL}}(2,\mbb{C})$. If $S/H$ is
K\"ahler, then $H$ is algebraic.
\end{lem}

\begin{proof}
Since every Lie subalgebra of $\lie{s}=\lie{sl}(2,\mbb{C})$ is conjugate to
$\{0\}$, to $\mbb{C}\left(\begin{smallmatrix}0&1\\0&0\end{smallmatrix}\right)$,
to $\mbb{C}\left(\begin{smallmatrix}1&0\\0&-1\end{smallmatrix}\right)$, to a
Borel subalgebra $\lie{b}$, or to $\lie{s}$, we conclude that the identity
component $H^0$ is automatically algebraic. Therefore it suffices to show that
$H$ has only finitely many connected components since then $H$ is the
finite union of translates of $H^0$ which is algebraic.

For $H^0=S$ this is trivial. Since the normalizer of a Borel subgroup $B$ of
$S$ coincides with $B$, we see that $H^0=B$ implies $H=B$, hence that $H$ is
algebraic in this case. If $H^0$ is a maximal algebraic torus in $S$, then its
normalizer in $S$ has two connected components, thus $H$ has at most two
connected components as well.

Suppose that $H^0$ is unipotent. Then its normalizer is a Borel subgroup. If
$H$ has infinitely many connected components, we find an element $h\in
H\setminus H^0$ which generates a closed infinite subgroup $\Gamma$ of $S$.
Then $S/(\Gamma H^0)$ is K\"ahler (for it covers $S/H$), and we conclude from
Lemma~\ref{Lem:twistedproduct} that $(S/H^0)\times_\Gamma \mbb{C}^*$ is a
Hamiltonian $S$--manifold where $\Gamma$ acts on $\mbb{C}^*$ by $\gamma^m\cdot
z:=e^{im}z$. As above this contradicts Theorem~\ref{Thm:locallyclosed}.

Since the case that $H^0$ is trivial, i.\,e.\ that $H$ is discrete, has already
been treated, the proof is finished.
\end{proof}

Now suppose that $h=h_{\sf s}h_{\sf u}$ with $h_{\sf s}\not=e$ and $h_{\sf u}
\not=e$. In this case there is a simple three dimensional closed complex
subgroup $S$ of the centralizer of $\mathcal{A}(h_{\sf s})$ which contains
$\mathcal{A}(h_{\sf u})$. Then $\mathcal{A}(h)\subset S\mathcal{A}(h_{\sf s})$
and a finite covering of $S\mathcal{A}(h_{\sf s})$ is isomorphic to
${\rm{SL}}(2,\mbb{C})\times(\mbb{C}^*)^l$. (We may suppose that
$\mathcal{A}(h_{\sf s})$ has positive dimension, if not, we are essentially in
the previous case.) Moreover, there is a closed complex subgroup $\wt{H}$ of
$\wt{G}={\rm{SL}}(2,\mbb{C})\times(\mbb{C}^*)^l$ containing the element
$h=\bigl(\left(\begin{smallmatrix}1&1\\0&1\end{smallmatrix}\right),
(e^{a_1},\dotsc,e^{a_l})\bigr)$ such that $\wt{G}/\wt{H}$ is Hamiltonian. We
must show that $\mathcal{A}(h)=\left(\begin{smallmatrix}1&\mbb{C}\\0&1
\end{smallmatrix}\right)\times(\mbb{C}^*)^l$ is contained in $\wt{H}$.

In order to simplify the notation we will continue to write $G$ and $H$ instead
of $\wt{G}$ and $\wt{H}$. The following observation is central to our argument.

\begin{lem}
We may assume without loss of generality that $H\cap(\mbb{C}^*)^l=\{e\}$.
\end{lem}

\begin{proof}
For this note that the action of $(\mbb{C}^*)^l$ on $G/H$ is Hamiltonian. This
implies that $H\cap(\mbb{C}^*)^l$ is a central subtorus $T$ of $G$.
Consequently, $G/H\cong(G/T)/(H/T)$. If $H/T$ is algebraic in $G/T$, then $H$
is algebraic in $G$.
\end{proof}

Let $p_1$ and $p_2$ denote the projections of $G={\rm{SL}}(2,\mbb{C})\times
(\mbb{C}^*)^l$ onto ${\rm{SL}}(2,\mbb{C})$ and $(\mbb{C}^*)^l$, respectively.

\begin{lem}
The map $p_1\colon G\to{\rm{SL}}(2,\mbb{C})$ maps the group $H$ isomorphically
onto a closed complex subgroup of ${\rm{SL}}(2,\mbb{C})$.
\end{lem}

\begin{proof}
We show first that $p_1(H)$ is closed in ${\rm{SL}}(2,\mbb{C})$. For this note
that $G/H$ is a Hamiltonian $(\mbb{C}^*)^l$--manifold. By
Theorem~\ref{Thm:locallyclosed} all $(\mbb{C}^*)^l$--orbits are locally
closed in $G/H$. Since $(\mbb{C}^*)^l$ is the center of $G$, we have
$(\mbb{C}^*)^l\cdot(gH)=g\cdot\bigl((\mbb{C}^*)^l\cdot eH\bigr)$. Hence all
$(\mbb{C}^*)^l$--orbits have the same dimension. This implies that all
$(\mbb{C}^*)^l$--orbits are closed in $G/H$. Consequently, $(\mbb{C}^*)^lH$ is
closed in $G$ which shows that $p_1(H)={\rm{SL}}(2,\mbb{C})\cap(\mbb{C}^*)^lH$
is closed.

Since the restriction of $p_1$ to the closed subgroup $H$ of $G$ is a
surjective holomorphic homomorphism onto $p_1(H)$ with kernel $H\cap
(\mbb{C}^*)^l=\{e\}$, the claim follows.
\end{proof}

If $p_1(H)={\rm{SL}}(2,\mbb{C})$, then $p_2\colon H\cong{\rm{SL}}(2,\mbb{C})\to
(\mbb{C}^*)^l$ must be trivial. But this contradicts the fact that $(e^{a_1},
\dotsc,e^{a_l})$ is contained in $p_2(H)$. Therefore $p_1(H)$ must be a proper
closed subgroup of ${\rm{SL}}(2,\mbb{C})$ which contains the element
$\left(\begin{smallmatrix}1&1\\0&1\end{smallmatrix}\right)$. In particular, we
conclude that $H^0$ is solvable.

There are essentially three possibilities. The image $p_1(H)$ is a Borel
subgroup of ${\rm{SL}}(2,\mbb{C})$ (which implies that $H$ is a connected
two-dimensional non-Abelian subgroup of $G$), or $p_1(H)^0=\left(
\begin{smallmatrix}1&\mbb{C}\\0&1\end{smallmatrix}\right)$, or $p_1(H)$ is
discrete containing $\left(\begin{smallmatrix}1&\mbb{Z}\\0&1\end{smallmatrix}
\right)$. If $p_1(H)$ is discrete, then $H$ is discrete. We have already shown
that $H$ is finite in this case, hence algebraic.

\begin{rem}
Suppose that $p_1(H)$ is the Borel subgroup of upper triangular matrices in
${\rm{SL}}(2,\mbb{C})$. The map $p_2|_H\colon H\to p_2(H)$ is a surjective
homomorphism with kernel $H'=H\cap{\rm{SL}}(2,\mbb{C})$. Thus we have
$H\cap{\rm{SL}}(2,\mbb{C})=\left(\begin{smallmatrix}1&\mbb{C}\\0&1
\end{smallmatrix}\right)$.
\end{rem}

Suppose that $p_1(H)$ is one-dimensional or a Borel subgroup. We know that
$p_2(H)$ is a closed complex subgroup of $(\mbb{C}^*)^l$ containing
$(e^{a_1},\dotsc,e^{a_l})$. Since $\dim p_2(H)=1$, we conclude that
$p_2(H)^0=\bigl\{(e^{ta_1},\dotsc,e^{ta_l});\ t\in\mbb{C}\bigr\}$.
Lemma~\ref{Lem:cover} implies that if $G/H$ is Hamiltonian, then the same holds
for $G/H^0$. The possibilities for $H^0$ are 
$H^0=\left\{\left(\left(\begin{smallmatrix}e^{ta_0}&s\\0&e^{-ta_0}
\end{smallmatrix}\right),(e^{ta_1},\dotsc,e^{ta_l}) \right);\
t,s\in\mbb{C}\right\}$ (if $p_1(H)$ is a Borel subgroup) or
$H^0=\left\{\left(\left(\begin{smallmatrix}
1&t\\0&1\end{smallmatrix}\right),(e^{ta_1},\dotsc,e^{ta_l})\right);\
t\in\mbb{C}\right\}$. We have to show that in both cases $G/H^0$ is not
Hamiltonian.

Let us first consider the case that
\begin{equation*}
H=\left\{\left(\begin{pmatrix}1&t\\0&1\end{pmatrix},(e^{ta_1},\dotsc,e^{ta_l})
\right);\ t\in\mbb{C}\right\}.
\end{equation*}
Let $T=(\mbb{C}^*)^{l-1}\times\{1\}\subset(\mbb{C}^*)^l$ and $\wt{G}:=
{\rm{SL}}(2,\mbb{C})\times T$. Then we have $T\cap p_2(H)=:\Gamma\cong\mbb{Z}$
and $(\mbb{C}^*)^l/p_2(H)\cong T/\Gamma$. Moreover, $G/H$ is a Hamiltonian
$\wt{G}$--manifold and we have $G/H\cong\wt{G}/\wt{H}$, where $\wt{H}=\wt{G}
\cap H\cong\mbb{Z}$. This contradicts our result in the discrete case. Hence,
$G/H$ cannot be Hamiltonian.

Finally, suppose that
\begin{equation*}
H=\left\{\left(\begin{pmatrix}e^{ta_0}&s\\0&e^{-ta_0}\end{pmatrix},
(e^{ta_1},\dotsc,e^{ta_l})\right);\ t,s\in\mbb{C}\right\}.
\end{equation*}
Again we consider $T=(\mbb{C}^*)^{l-1}\times\{1\}$ and $\wt{G}={\rm{SL}}(2,
\mbb{C})\times T$. We have $\wt{H}=\wt{G}\cap H=\Gamma\ltimes\left(
\begin{smallmatrix}1&\mbb{C}\\0&1\end{smallmatrix}\right)$ where $\Gamma\cong
\mbb{Z}$. As in the discrete case we let $\Gamma$ act on $\mbb{C}^*$ by
$\gamma^m\cdot z:=e^{im}z$ and consider the twisted product $(\wt{G}/\wt{H}^0)
\times_\Gamma\mbb{C}^*$. If $G/H$ is Hamiltonian, then the same holds for
$\wt{G}/\wt{H}$ and thus for $\wt{G}/\wt{H}^0$. Then
$(\wt{G}/\wt{H}^0)\times_\Gamma\mbb{C}^*$ is Hamiltonian by
Lemma~\ref{Lem:twistedproduct}. Since the $\Gamma$--orbits in $\mbb{C}^*$ are
not locally closed, this contradicts Theorem~\ref{Thm:locallyclosed} and we
conclude that $G/H$ is not Hamiltonian.

Summarizing our discussion in this section we proved the following.

\begin{thm}
Let $G$ be a connected complex reductive  group and let $H$ be a closed complex
subgroup. If $X=G/H$ is a Hamiltonian $G$--manifold, then $H$ is an algebraic
subgroup of $G$.
\end{thm}

In particular we obtain the following result which was originally proved
in~\cite{BerO} and~\cite{Bert}.

\begin{cor}\label{Cor:SemisimpleCase}
Let $S$ be a connected semisimple complex Lie group and let $H$ be a closed
complex subgroup of $S$. If $S/H$ admits a K\"ahler form, then $H$ is an
algebraic subgroup of $S$.
\end{cor}

\section{Homogeneous K\"ahler manifolds}

Let $G=K^\mbb{C}$ be a connected complex reductive  Lie group. In this section
we characterize those closed complex subgroups $H$ of $G$ for which $X=G/H$
admits a K\"ahler form.

According to~\cite{Bo}, Corollary~I.2.3, the commutator group $S:=G'$ is a
connected algebraic subgroup of $G$ and, since $G$ is reductive, $S$ is
semisimple. Let $Z:=\mathcal{Z}(G)^0\cong(\mbb{C}^*)^k$. Then $G=SZ$ and
$S\cap Z$ is finite.

\begin{thm}\label{Thm:Kaehler}
Let $G$ be a reductive complex Lie group and $H\subset G$ a closed complex
subgroup. Then the manifold $X=G/H$ admits a K\"ahler form if and only if $S
\cap H\subset S$ is algebraic and $SH$ is closed in $G$.
\end{thm}

\begin{proof}
After replacing $G$ by a finite cover we may assume that $G=S\times Z$. Suppose
first that $G/H$ is K\"ahler. Then $S/(S\cap H)$ is also K\"ahler and hence
$S\cap H\subset S$ is algebraic by Corollary~\ref{Cor:SemisimpleCase}. By
Theorem~\ref{Thm:locallyclosed} all $S$--orbits in $G/H$ are open in their
closures and their boundaries only contain orbits of strictly smaller
dimension. In view of the reductive group structure of $G$ the $S$--orbits in
$X$ all have the same dimension. This implies that every $S$--orbit in $G/H$ is
closed. Consequently, $SH$ is closed in $G$.

Now suppose that $S\cap H\subset S$ is algebraic and that $SH$ is closed in
$G$. Although it is not used in the proof we remark that we may assume that $H$
is solvable, since otherwise one can divide by the (ineffective) semisimple
factor of $H$. Consider the fibration
\begin{equation*}
X=G/H\to G/SH.
\end{equation*}
The base $G/SH$ is an Abelian complex Lie group and the fiber is $SH/H=S/(S\cap
H)$. There is a subgroup $(\mbb{C}^*)^l\cong Z_1\subset Z\cong(\mbb{C}^*)^k$
such that $G_1:=S\times Z_1\subset G$ acts transitively on $X$ and $Z_1\cap SH$
is discrete. With $H_1:=H\cap G_1$ we have that $X=G_1/H_1$, that $SH_1$ is
closed in $G_1$ and that the base of the fibration 
\begin{equation*}
X=G_1/H_1\to G_1/SH_1=Z_1/(Z_1\cap SH_1)
\end{equation*}
is a discrete quotient of $Z_1$. Furthermore $S\cap H=S\cap H_1$ is an
algebraic subgroup of $S$.

Let $\Gamma_1:=Z_1\cap SH_1$ and $\Gamma_2\subset Z_1$ be a discrete subgroup
such that $\Gamma_1\cap\Gamma_2=\{e\}$ and such that $\Gamma:=\Gamma_1+
\Gamma_2$ is a discrete cocompact subgroup of $Z_1$. Since $H_1'$ is contained
in $S\cap H_1$, we can define the closed complex subgroup $H_2\subset G_1$ to
be the group generated by $H_1$ and $\Gamma_2(S\cap H_1)$. One still has that
$H_2\cap S=H_1\cap S=H\cap S$ is algebraic in $S$ and that $SH_2$ is closed in
$G_1$. Hence $X=G_1/H_1\to G_1/H_2$ is a covering map and one sees that
in order to finish the proof it is sufficient to show that the base $G_1/H_2$
admits a K\"ahler form.

So we may drop the indices and have to prove that a reductive quotient $X=G/H$
of $G=S\times(\mbb{C}^*)^k=S\times Z$  by a closed complex solvable subgroup
$H$ with $S\cap H\subset S$ algebraic, $SH\subset G$ closed and $G/SH$ a
compact torus, is K\"ahler.

Let ${\rm p}_1\colon G\to S$ be the projection onto $S$. Note that the
algebraic Zariski closure $\ol{H}$ of $H$ in $G$ is the product $\wh{H}\times
(\mbb{C}^*)^k$, where $\wh{H}$ is the Zariski closure of the projection ${\rm
p}_1(H)$ of $H$ in $S$. The commutator group $H'$ of $H$ is also the commutator
group of $\ol{H}$ and of $\wh{H}$  and is contained in $H\cap S$. Therefore one
gets a natural algebraic right action of $\ol{H}$ on the homogeneous manifold
$Y:=S/(S\cap H)$ given by
\begin{equation*}
\ol{h}\bigl(s(S\cap H)\bigr):=s\bigl({\rm p}_1(\ol{h})\bigr)^{-1}(S\cap H).
\leqno{(\star)}
\end{equation*}
As a consequence we can equivariantly compactify the $S\ol{H}$--manifold $Y$ to
an almost homogeneous projective $S\ol{H}$--manifold $\ol{Y}$, see \cite{Mass},
Proposition~3.1, and \cite{Koll}, Proposition~3.9.1. Since $X=G/H$ is
realizable as a quotient of the manifold $S/(S\cap H)\times(\mbb{C}^*)^k$ by
the natural action of $H/(S\cap H)$, where the $S$--factor of the  action is
given by $(\star)$, we see that $X$ is an open set in a holomorphic fiber
bundle $\ol{X}$ with a compact torus as base and the simply connected
projective manifold $\ol{Y}$ as fiber. Finally we apply Blanchard's theorem,
see~\cite{Bl}, p.~192, to get a K\"ahler form on $\ol{X}$ and then, by
restriction, on $X$ also. The theorem is proved.
\end{proof}

\end{document}